\documentclass[final,1p,times]{elsarticle}%
\usepackage{lineno,hyperref}
\usepackage{graphicx}
\usepackage{amsmath}
\usepackage{amsfonts}
\usepackage{amssymb}
\usepackage{graphicx}
\usepackage{multirow}
\usepackage{enumitem}
\usepackage{color}%

\usepackage{amsmath,multirow}
\usepackage{amssymb}
\usepackage{amscd}
\usepackage{amsfonts}
\usepackage{enumerate}
\usepackage{mathrsfs}
\usepackage{xypic}
\usepackage{stmaryrd}
\usepackage{graphicx}
\usepackage{fancybox}
\usepackage{portland}
\usepackage{color}
\usepackage{amsmath}
\usepackage{exscale}
\usepackage{enumitem,array}%
\usepackage{amssymb}
\newtheorem{mytheo}{Theorem}[section]
\newtheorem{mydef}[mytheo]{Definition}

\newtheorem{rem}{Remark}

\journal{}

\begin{document}

\begin{frontmatter}

\title{Efficient implementation of RKN-type Fourier collocation methods for second-order differential equations\tnoteref{mytitlenote}} \tnotetext[mytitlenote]{This paper
was supported by  National Natural Science Foundation of China
(Grant No. 11401333, 11171178, 11571302), by  Natural Science
Foundation of Shandong Province (Grant No. ZR2014AQ003), by
Foundation of Scientific Research Project of Shandong Universities
(Grant No. J14LI04), and by  China Postdoctoral Science Foundation
(Grant No. 2015M580578). }

\author[Wang]{Bin~Wang\corref{cor1}}
\author[Wang]{Fanwei Meng}
\author[Fang]{Yonglei Fang}
\address[Wang]{School of Mathematical Sciences, Qufu Normal University, Qufu, Shandong 273165, P.R.China}
\address[Fang]{School of Mathematics  and Statistics, Zaozhuang University,  Zaozhuang, Shandong  277160, P.R. China}

 \ead{wangbinmaths@gmail.com (Bin Wang), fwmeng@mail.qfnu.edu.cn (Fanwei Meng), ylfangmath@163.com}
\cortext[cor1]{Corresponding author.}

\begin{abstract}
In this paper we discuss the efficient implementation of RKN-type
Fourier collocation methods, which are  used when solving
second-order differential equations. The proposed implementation
relies on an alternative   formulation of the methods and the
blended formulation. The features and effectiveness of the
implementation are confirmed by the performance of the methods on
two numerical tests.
\end{abstract}


\end{frontmatter}

\section{Introduction}
 The efficient numerical solution of implicit methods  for
differential equations has been the subject of many investigations
in the last decades. This paper is devoted to dealing with the
numerical solution of  second-order differential equations, namely
problems in the form
\begin{equation}
q^{\prime\prime}(t)=f(q(t)),  \qquad
q(t_0)=q_0,\ \ q'(t_0)=q_0',\qquad t\in[t_0,t_{\mathrm{end}}],\label{common prob}%
\end{equation}
where $q(t): \mathbb{R}\rightarrow \mathbb{R}^{d}$ and  $f(q):
\mathbb{R}^{d}\rightarrow \mathbb{R}^{d}$ is an analytic function.
The solution of this system and  its derivative satisfy the
following variation-of-constants formula (\cite{wu2013-book}) with
the stepsize $h$
\begin{equation}\label{variation-of-constants formula}
\begin{aligned} &q(t+\mu h)=q(t)+\mu
hq'(t)+h^2\int_0^\mu
(\mu-x)f(q(t+hx))\,\mathrm{d}x,\\
&q'(t+\mu h)=q'(t)+h\int_0^\mu f(q(t+hx))\,\mathrm{d}x.
\end{aligned}
\end{equation}

Numerical methods of  the second-order system \eqref{common prob}
have been studied by many researchers in the last decades (see, e.g.
\cite{Butcher-book,franco2002,Skeel-1998,hairer2006,Simos-2012,Simos2002,ANM-2015,wu2013-ANM,wu-2012-BIT,Tang-2004}),
and Runge--Kutta--Nystr\"{o}m (RKN) methods   are one of well-known
methods for solving this system.
 In \cite{wang-2016},  the authors  took advantage of shifted
Legendre polynomials to obtain a local Fourier expansion of the
considered system and derived a  kind of collocation methods
(trigonometric Fourier collocation methods). The analysis given in
\cite{wang-2016} also presents a new collocation methods (RKN-type
Fourier collocation methods) for solving the second-order system
\eqref{common prob}. This kind of collocation methods is implicit
and an
 iterative procedure is required.
 On the other hand, efficient implementation of   implicit
 methods has been investigated by many researchers  in recent years and we refer  to
 \cite{Brugnano2013,Brugnano2014,Brugnano2015,Brugnano2011,Brugnano2002,Brugnano2004,Brugnano2007}
 for some examples on this topic.   Motivated by these publications,
  here we present  an efficient implementation of RKN-type Fourier collocation  methods, by proposing and analysing
an iterative procedure based on the particular structure of the
methods.

With this premise, this paper is organized as follows.  Section
\ref{sec:RKN-type  Fourier collocation method}   describes the
derivation of RKN-type Fourier collocation methods and the structure
of the discrete problem generated by the methods. In Section
\ref{sec:Implementation} we propose and analyze an efficient
implementation of RKN-type Fourier collocation methods. Section
\ref{sec:Numerical tests} reports two numerical tests to show the
features and effectiveness of the methods. The last section contains
a few conclusions.

\section{RKN-type  Fourier collocation
methods} \label{sec:RKN-type  Fourier collocation method} RKN-type
Fourier collocation methods are given in \cite{wang-2016} as a
by-product of trigonometric Fourier collocation methods for
second-order differential equations.
 We now recall the   derivation of the
methods, and derive the most efficient formulation of the generated
discrete problems.

 Let  us   consider the restriction of problem
\eqref{common prob} to the interval $[t_0,t_0+h]$, with the
right-hand side expanded along the shifted Legendre polynomials
$\{\widehat{P}_j\}_{j=0}^{\infty}$ over the interval $[0,1]$, scaled
in order to be orthonormal. Then we rewrite the function $f(q)$ in
\eqref{common prob} as
\begin{equation*}f(q(t_0+\xi h))=\sum\limits_{j=0}^
{\infty}\widehat{P}_j(\xi
)\int_{0}^{1}\widehat{P}_j(\tau)f(q(t_0+\tau h))d\tau,\ \ \xi
\in[0,1].
\end{equation*}
Thence the problem \eqref{common prob} to the interval $[t_0,t_0+h]$
is rewritten as
\begin{equation*}\begin{aligned}
q^{\prime\prime}(t_0+\xi h) =&\sum\limits_{j=0}^
{\infty}\widehat{P}_j(\xi
)\int_{0}^{1}\widehat{P}_j(\tau)f(q(t_0+\tau h))d\tau,\\&  \ q(t_0)=q_0,\ \ q'(t_0)=q_0'.\label{expanded prob}%
\end{aligned}\end{equation*}
Truncating the series at the right-hand side gives the following
approximate problem
\begin{equation*}\begin{aligned}
u^{\prime\prime}(t_0+\xi h)=&\sum\limits_{j=0}^
{r-1}\widehat{P}_j(\xi
)\int_{0}^{1}\widehat{P}_j(\tau)f(u(t_0+\tau h))d\tau,\\& \ u(t_0)=q_0,\ \ u'(t_0)=q_0'.\label{approximate prob}%
\end{aligned}\end{equation*}

 Introduce a quadrature formula based at $k\ (k\geq r)$ abscissae
$0\leq c_1\leq\ldots\leq c_k\leq1$ to deal with
$\int_{0}^{1}\widehat{P}_j(\tau)f(u(t_0+\tau h))d\tau$. Thus we
obtain an approximation of the form
\begin{equation*}
\begin{aligned} \int_{0}^{1}\widehat{P}_j(\tau)f(u(t_0+\tau h))d\tau\approx \sum\limits_{l=1}^
{k}b_l\widehat{P}_j(c_l)f(u(t_0+c_l h)),\ \ j=0,1,\ldots, r-1,
\end{aligned}
\label{Computate gamma}%
\end{equation*} where $b_l$ with $l=1,2,\ldots,k$ are the
quadrature weights. It is natural to consider the following discrete
problem  as an approximation of \eqref{common prob}
\begin{equation}
\begin{aligned}
&v''(t_0+\xi h)=\sum\limits_{j=0}^ {r-1}\widehat{P}_j(\xi
)\sum\limits_{l=1}^
{k}b_l\widehat{P}_j(c_l)f(v(t_0+c_l h)),\ \ \  v(t_0)=q_{0},\ v'(t_0)=q_0',\\
\end{aligned}
\label{discrete quadrature truncating H-s}%
\end{equation}
which can be solved by the variation-of-constants formula
\eqref{variation-of-constants formula} in the form:
\begin{equation*}
\begin{aligned}
v_i=&q_0+c_ihq'_0+ h^2\sum\limits_{j=0}^
{r-1}\int_{0}^{c_i}\widehat{P}_j(x)(c_i-x)dx
\sum\limits_{l=1}^ {k}b_l\widehat{P}_j(c_l)f(v_l),\\
i=&1,2,\ldots,k,\\
\end{aligned}
\end{equation*}
where $v_i= v(t_0+c_ih).$

The approximation to $q(t_0+ h),\ q'(t_0+ h)$ is then given by
$q_1=v(t_0+ h),\ q'_1=v'(t_0+ h),$ which can be obtained by applying
the variation-of-constants formula \eqref{variation-of-constants
formula} to \eqref{discrete quadrature truncating H-s} as follows:
\begin{equation*}
\begin{aligned}
q_1=&q_0+hq'_0+h^2\sum\limits_{j=0}^
{r-1}\int_{0}^1\widehat{P}_j(x)(1-x)dx
\sum\limits_{l=1}^ {k}b_l\widehat{P}_j(c_l)f(v_l),\\
q'_1=&q'_0 +h\int_{0}^1\widehat{P}_j(x)dx \sum\limits_{l=1}^ {k}b_l
f(v_l).
\end{aligned}
\end{equation*}
We compute
\begin{equation*}
\begin{aligned}
&\int_{0}^1\widehat{P}_j(x)(1-x)dx =\left\{
\begin{aligned}
 &\frac{1}{2},\ j=0,\\
 &\frac{-1}{2\sqrt{3}}, \ j=1,\\
  &0, \ j\geq2,\\\end{aligned}\right. \qquad
  \int_{0}^1\widehat{P}_j(x)dx
=\left\{
\begin{aligned}
 &1,\ j=0,\\
 &0, \ j\geq1.\end{aligned}\right.\\
\end{aligned}
\end{equation*}
Then we have the following result for $r\geq2:$
\begin{equation*}
\begin{aligned}
q_1=&q_0+hq'_0+ h^2 \sum\limits_{l=1}^
{k}(1-c_l)b_lf(v_l),\\
q'_1=&q'_0 +h \sum\limits_{l=1}^ {k}b_l f(v_l).
\end{aligned}
\end{equation*}

 We are now in a position to present  RKN-type  Fourier
collocation methods for the second-order ODEs \eqref{common prob}.
 \begin{mydef}
\label{numerical RKN method} (\cite{wang-2016})  A $k$-stage
RKN-type  Fourier collocation method for integrating the system
\eqref{common prob} is defined as
\begin{equation}
\begin{aligned} v_i&=q_0+c_ihq'_0+
h^2\sum\limits_{j=0}^ {r-1}\int_{0}^{c_i}\widehat{P}_j(x)(c_i-x)dx
\sum\limits_{l=1}^ {k}b_l\widehat{P}_j(c_l)f(v_l),\
i=1,2,\ldots,k,\\
q_1&=q_0+hq'_0+ h^2 \sum\limits_{l=1}^
{k}(1-c_l)b_lf(v_l),\\
q'_1&=q'_0 +h \sum\limits_{l=1}^ {k}b_l f(v_l),
\end{aligned}
\label{methods0}%
\end{equation}
where $h$ is the stepsize,  $r$ is an integer with the requirement
$2\leq r\leq k$, $(c_l,b_l)$ with $l=1,2,\ldots,k$ are the node
points and the quadrature weights of a quadrature formula,
respectively.
\end{mydef}
\begin{rem} It is remarked that the properties of the methods including the convergence,
stability, and  the degree of accuracy    in preserving the
solution, the quadratic invariant
  and the Hamiltonian   have been studied in \cite{wang-2016}. \end{rem}

 It is noted that
the  method \eqref{methods0} is the subclass of $k$-stage RKN
methods with the following Butcher tableau:
\begin{equation} \label{Butcher of RKN methods}
\begin{tabular}
[c]{l}%
\\
\\[2mm]%
\begin{tabular}
[c]{c|c}%
$c$ & $\bar{A}=(\bar{a}_{lm})_{k\times k}$\\\hline &
$\raisebox{-1.3ex}[0pt]{$\bar{b}^{\intercal}$}$\\\hline
& $\raisebox{-1.3ex}[0.5pt]{$b^{\intercal}$}$%
\end{tabular}
$\ \quad=$ $\ $%
\end{tabular}%
\begin{tabular}
[c]{c|ccc}%
$c_{1}$  \\
$\vdots$ &    $ \bar{a}_{lm} =b_m\sum\limits_{j=0}^
{r-1}\int_{0}^{c_l}\widehat{P}_j(x)(c_l-x)dx
 \widehat{P}_j(c_m)$\\
$c_{k}$ \\\hline & $\raisebox{-1.3ex}[0.5pt]{$(1-c_1)b_1\ \ \cdots\
\ (1-c_k)b_k$}$
\\\hline & $\raisebox{-1.3ex}[1.0pt]{$b_1\ \ \cdots\ \ b_k$}$
\end{tabular}
\end{equation}
  It is convenient to
express the  methods in block-matrix notation
\begin{equation}\label{block-matrix of methods}
\begin{aligned}
&v=u\otimes q_0+hc\otimes q'_0+h^2\bar{A}\otimes I_df(v),\\
&q_1=q_0+hq'_0+h^2\left(\bar b^{\intercal}\otimes I_d\right)f(v),\\
&q'_1=q'_0 +h\left(b^{\intercal}\otimes I_d\right)f(v),
\end{aligned}\end{equation} where $I_d$ is the $d\times d$ identity
matrix, $\otimes$ is the Kronecker product,
$u=\left(1,\cdots,1\right)^{\intercal}$,
$v=\left(v_1^{\intercal},\cdots,v_k^{\intercal}\right)^{\intercal}$
and
$f(v)=\big(f(v_1)^{\intercal},\cdots,f(v_k)^{\intercal}\big)^{\intercal}$.

 It can be observed that usually the method \eqref{methods0}
constitutes of a system of implicit equations for the determination
of $v_i\ (i=1,2,\cdots,k)$ and it requires an iterative  procedure.
It is noted that usually $k\geq r$ in \eqref{methods0} and the
iterative computation can be reduced to solve an implicit system
with $r$ equations, which is quite important for practical
computations especially when $r<k$. The effective implementation of
RKN-type Fourier collocation method will be discussed in next
section.


\section{Implementation of the methods} \label{sec:Implementation}
In this section, we propose and analyze the efficient implementation
of RKN-type Fourier collocation methods.
\subsection{Fundamental and silent stages} \label{subsec:Fundamental and silent stages}
By introducing the matrices
$\Omega_k=\mathrm{diag}(b_1,b_2,\cdots,b_k),$
\begin{equation}
 \mathcal{L}_{k,r}=\Bigg(\int_{0}^{c_i}\widehat{P}_j(x)(c_i-x)dx\Bigg)_{\begin{subarray}{l}
                                                             i=1,2,\cdots,k \\
                                                             j=0,1,\cdots,r-1
                                                             \end{subarray}}
                                                            \in\mathbb{R}^{k\times r},
\label{Lkr}%
\end{equation}
and
\begin{equation}
 \mathcal{P}_{k,r}=\Big(\widehat{P}_j(c_i)\Big)_{\begin{subarray}{l}
                                                            i=1,2,\cdots,k \\
                                                          j=0,1,\cdots,r-1
                                                            \end{subarray}}
                                                            \in\mathbb{R}^{k\times
                                                            r},\label{Pkr}%
\end{equation}
the   matrix $\bar{A}$ of RKN-type  Fourier collocation method given
in \eqref{Butcher of RKN methods} can be recast as
\begin{equation}
\bar{A}= \mathcal{L}_{k,r} \mathcal{P}^{\intercal}_{k,r}\Omega_k.\label{barA}%
\end{equation}
\begin{mytheo}\label{LKR  thm}
The matrix $\bar{A}$ in \eqref{barA} can be written as
\begin{equation}
\bar{A}= \mathcal{P}_{k,r+2} \left(
                                                                           \begin{array}{c}
                                                                             X_{r,r} \\
                                                                             \xi_{r-1}\xi_{r}e^{\intercal}_{r-1} \\
                                                                             \xi_{r}\xi_{r+1}e^{\intercal}_{r}

                                                                           \end{array}
                                                                         \right) \mathcal{P}^{\intercal}_{k,r}\Omega_k,\label{new barA}%
\end{equation}
where  $$\xi_{m}=\frac{1}{2\sqrt{4m^2-1}},$$   \begin{equation}\label{XRR}\begin{aligned} &X_{r,r}= \\
&\left(
                                            \begin{array}{cccccccc}
                                            \frac{1}{4}-\xi_1^2 & -\frac{\xi_1}{2} & \xi_1\xi_2   &   &  & &  \\
                                            \frac{\xi_1}{2} & -\xi_1^2-\xi_2^2 & 0 & \xi_2\xi_3    &  &  &   \\
                                            \xi_1\xi_2 & 0 &  -\xi_2^2-\xi_3^2  & 0 &  \ddots     &  &  \\
                                             & \xi_2\xi_3 & 0 & \ddots   & \ddots &   \ddots   &      \\
                                            &  &  \ddots & \ddots  &  \ddots & 0  &  \xi_{r-2}\xi_{r-1} \\
                                              &  &  &      \ddots& 0  & -\xi_{r-2}^2-\xi_{r-1}^2  & 0 \\
                                            &  & &     & \xi_{r-2}\xi_{r-1}  &  0  & -\xi_{r-1}^2-\xi_r^2\\
                                            \end{array}
                                          \right),\end{aligned}
\end{equation} and
$e_{r-1},\ e_{r}$ are the unit coordinate vectors of length $r-1$
and $r$, respectively.
\end{mytheo}
\textbf{Proof.}  According to \cite{hairer-stiff}, the shifted
Legendre polynomials $\{\widehat{P}_j\}_{j=0}^{\infty}$ satisfy the
 following  integration formulae
\begin{equation}
\begin{aligned}\label{L polynomials P}%
&\int_{0}^{x}\widehat{P}_0( t)dt=\xi_1\widehat{P}_1( x)+\frac{1}{2}\widehat{P}_0( x),\\
 &\int_{0}^{x}\widehat{P}_m( t)dt=\xi_{m+1}\widehat{P}_{m+1}( x)-\xi_{m}\widehat{P}_{m-1}( x),\ m=1,2,\ldots.\end{aligned}
\end{equation}
 This yields that
\begin{equation*}
\begin{aligned}
&\widehat{P}_0(x)=\xi_1\widehat{P}'_1(x)+\frac{1}{2}\widehat{P}_0'( x),\\
 &\widehat{P}_m(x)=\xi_{m+1}\widehat{P}'_{m+1}( x)-\xi_{m}\widehat{P}'_{m-1}( x),\ m=1,2,\ldots.\end{aligned}
\end{equation*}
We compute
\begin{equation*}
\begin{aligned}
&\int_{0}^{c_i}\widehat{P}_0(x)(c_i-x)dx=c_i\int_{0}^{c_i}\widehat{P}_0(x)dx-\int_{0}^{c_i}x\widehat{P}_0(x)dx\\
=&c_i\int_{0}^{c_i}\widehat{P}_0(x)dx-\int_{0}^{c_i}x\Big(\xi_1\widehat{P}'_1(x)+\frac{1}{2}\widehat{P}_0'( x)\Big)dx\\
=&c_i\int_{0}^{c_i}\widehat{P}_0(x)dx-\xi_{1}\int_{0}^{c_i}xd\widehat{P}_{1}(
x)+\frac{1}{2}\int_{0}^{c_i}xd\widehat{P}_{0}( x)\\
=&c_i\int_{0}^{c_i}\widehat{P}_0(x)dx-\xi_{1}c_i\widehat{P}_{1}(
c_i)+\xi_{1}\int_{0}^{c_i}\widehat{P}_{1}(
x)dx+\frac{1}{2}c_i\widehat{P}_{0}(c_i)-\frac{1}{2}\int_{0}^{c_i}\widehat{P}_{0}(
x)dx,\end{aligned}
\end{equation*}and
\begin{equation*}
\begin{aligned}
&\int_{0}^{c_i}\widehat{P}_j(x)(c_i-x)dx=c_i\int_{0}^{c_i}\widehat{P}_j(x)dx-\int_{0}^{c_i}x\widehat{P}_j(x)dx\\
=&c_i\int_{0}^{c_i}\widehat{P}_j(x)dx-\int_{0}^{c_i}x\Big(\xi_{j+1}\widehat{P}'_{j+1}(
x)-\xi_{j}\widehat{P}'_{j-1}( x)\Big)dx\\
=&c_i\int_{0}^{c_i}\widehat{P}_j(x)dx-\xi_{j+1}\int_{0}^{c_i}xd\widehat{P}_{j+1}(
x)+\xi_{j}\int_{0}^{c_i}xd\widehat{P}_{j-1}( x)\\
=&c_i\int_{0}^{c_i}\widehat{P}_j(x)dx-\xi_{j+1}c_i\widehat{P}_{j+1}(
c_i)+\xi_{j+1}\int_{0}^{c_i}\widehat{P}_{j+1}(
x)dx+\xi_{j}c_i\widehat{P}_{j-1}(c_i)\\&-\xi_{j}\int_{0}^{c_i}\widehat{P}_{j-1}(
x)dx,\ \ \ \ m=1,2,\ldots.\end{aligned}
\end{equation*}
According to \eqref{L polynomials P}, we obtain
\begin{equation*}
\begin{aligned}\label{Lx polynomials P}%
\int_{0}^{c_i}\widehat{P}_0(x)(c_i-x)dx=&(\frac{1}{4}-\xi_1^2)\widehat{P}_0(c_i)+\frac{\xi_1}{2}\widehat{P}_1(c_i)+\xi_1\xi_2\widehat{P}_2(c_i),\\
\int_{0}^{c_i}\widehat{P}_1(x)(c_i-x)dx=&-\frac{\xi_1}{2}\widehat{P}_0(c_i)-(\xi_1^2+\xi_2^2)\widehat{P}_1(c_i)+\xi_2\xi_3\widehat{P}_2(c_i),\\
 \int_{0}^{c_i}\widehat{P}_j(x)(c_i-x)dx=&\xi_{j-1}\xi_j\widehat{P}_{j-2}(c_i)-(\xi_j^2+\xi_{j+1}^2)\widehat{P}_j(c_i)+\xi_{j+1}\xi_{j+2}\widehat{P}_{j+2}(c_i),\\
  &j=2,3,\ldots.\end{aligned}
\end{equation*}
It follows from these formulae that  the matrix $\mathcal{L}_{k,r}$
\eqref{Lkr} can be rewritten as
\begin{equation*}\label{Lkr P}%
\mathcal{L}_{k,r}=\mathcal{P}_{k,r+2}\hat{X}_{r+2,r}=\mathcal{P}_{k,r+2}
\left(
                                                                           \begin{array}{c}
                                                                             X_{r,r} \\
                                                                             \xi_{r-1}\xi_{r}e^{\intercal}_{r-1} \\
                                                                             \xi_{r}\xi_{r+1}e^{\intercal}_{r}

                                                                           \end{array}
                                                                         \right).
\end{equation*}
Then \eqref{new barA} holds true. \hfill  $\blacksquare$

We have the following result about the matrix $X_{r,r}$.
\begin{mytheo}\label{detX}
Let $S_r=\mathrm{det}(X_{r,r}).$ Then $\{S_r\}_{r=1}^{\infty}$
satisfy  the recursion
\begin{equation*}\begin{aligned} S_1&=\frac{1}{4}-\xi_1^2,\\
 S_{2n}&=-\xi_{2n}^2S_{2n-1}+\xi_{1}^4\xi_{3}^4\cdots\xi_{2n-1}^4,\\
  S_{2n+1}&=-\xi_{2n+1}^2S_{2n}+\frac{1}{4}\xi_{2}^4\xi_{4}^4\cdots\xi_{2n}^4,\ \ n=1,2,\ldots.\\
\end{aligned}
\end{equation*}
\end{mytheo}
\textbf{Proof.} The result easily follows from the Laplace
expansion, by considering that, from \eqref{XRR},
$S_1=\frac{1}{4}-\xi_1^2$ and
$S_2=\xi_1^2\xi_2^2-\frac{1}{4}\xi_2^2+\xi_1^4.$ \hfill
$\blacksquare$

From the above analysis, it follows  that a $k$-stage RKN-type
Fourier collocation method, with $k>r$, is defined by a Butcher
matrix of rank $r$. Therefore the discrete problem  can be recast in
a more convenient form, whose (block) size is $r$, rather than $k$.
For this purpose, let us partition the abscissae $c_i,\
i=1,2,\ldots,k$ into two sets: one with $r$ abscissae, the other
with the remaining $k-r$ ones. We choose them as the first $r$ ones
and the remaining $k-r$ ones, respectively, for the sake of
simplicity. According to \cite{Brugnano2014,Brugnano2011}, the
corresponding stages are called fundamental stages and silent
stages, respectively.
 The key
idea is now that the $k-r$  silent stages can be expressed as  a
linear combination of the $r$  fundamental ones.  To this end, let
us then partition the matrices   as follows:
\begin{equation*}
\mathcal{L}_{k,r} =\left(
                     \begin{array}{c}
                       \mathcal{L}^{(1)}_{k,r}\\
                       \mathcal{L}^{(2)}_{k,r}  \\
                     \end{array}
                   \right),\ \ \mathcal{P}_{k,r} =\left(
                     \begin{array}{c}
                       \mathcal{P}^{(1)}_{k,r} \\
                       \mathcal{P}^{(2)}_{k,r}  \\
                     \end{array}
                   \right),\ \ \Omega_k = \left(
                                            \begin{array}{cc}
                                             \Omega^{(1)}_{k}   &  \\
                                                &   \Omega^{(2)}_{k} \\
                                            \end{array}
                                          \right),
\end{equation*}
where $$\mathcal{L}^{(1)}_{k,r},\mathcal{P}^{(1)}_{k,r},
\Omega^{(1)}_{k} \in \mathbb{R}^{r\times r},\
\mathcal{L}^{(2)}_{k,r}, \mathcal{P}^{(2)}_{k,r}\in
\mathbb{R}^{(k-r)\times r},\   \Omega^{(2)}_{k}\in
\mathbb{R}^{(k-r)\times (k-r)}.$$
 Similarly, let us denote by $v^{(1)}$ the (block) vector, of
dimension $r$, containing the fundamental stages, and by $v^{(2)}$
the (block) vector, of dimension $k-r$, with the silent stages. One
then obtains the equations:
\begin{equation*}
\left(
  \begin{array}{c}
    v^{(1)} \\
    v^{(2)} \\
  \end{array}
\right)=u\otimes q_0+hc\otimes q'_0+h^2\left(
                     \begin{array}{c}
                       \mathcal{L}^{(1)}_{k,r} \\
                       \mathcal{L}^{(2)}_{k,r}  \\
                     \end{array}
                   \right)
\left(
                     \begin{array}{c}
                       \mathcal{P}^{(1)}_{k,r} \\
                       \mathcal{P}^{(2)}_{k,r}  \\
                     \end{array}
                   \right)^{\intercal}\left(
                                            \begin{array}{cc}
                                             \Omega^{(1)}_{k}   &  \\
                                                &   \Omega^{(2)}_{k} \\
                                            \end{array}
                                          \right)\otimes I_d\left(
                                                               \begin{array}{c}
                                                                 f(v^{(1)} ) \\
                                                                 f(v^{(2)} ) \\
                                                               \end{array}
                                                             \right),
\end{equation*}
which can be written as
\begin{equation}
\begin{aligned}
&v^{(1)}=u^{(1)}\otimes q_0+hc^{(1)}\otimes q'_0+h^2
\mathcal{L}^{(1)}_{k,r}\mathcal{P}^{\intercal}_{k,r}\Omega_k\otimes
I_d\left(
                                                               \begin{array}{c}
                                                                 f(v^{(1)} ) \\
                                                                 f(v^{(2)} ) \\
                                                               \end{array}
                                                             \right),\\
&v^{(2)}=u^{(2)}\otimes q_0+hc^{(2)}\otimes q'_0+h^2
\mathcal{L}^{(2)}_{k,r} \mathcal{P}^{\intercal}_{k,r}\Omega_k\otimes
I_d\left(
                                                               \begin{array}{c}
                                                                 f(v^{(1)} ) \\
                                                                 f(v^{(2)} ) \\
                                                               \end{array}
                                                             \right),\\
\end{aligned}\label{v1v2}%
\end{equation}
where $c^{(1)}=(c_1,\ldots,c_r)^{\intercal}$,
$c^{(2)}=(c_{r+1},\ldots,c_{k})^{\intercal}$, $u^{(1)}$ and
$u^{(2)}$ are the unit vectors of length $r$ and $k-r$,
respectively.

From the first formula of \eqref{v1v2}, it follows that
\begin{equation*}\mathcal{P}^{\intercal}_{k,r}\Omega_k\otimes
I_d\left(
                                                               \begin{array}{c}
                                                                 f(v^{(1)} ) \\
                                                                 f(v^{(2)} ) \\
                                                               \end{array}
                                                             \right)=(h^2
\mathcal{L}^{(1)}_{k,r})^{-1}\otimes I_d [v^{(1)} -u^{(1)}\otimes
q_0-hc^{(1)}\otimes q'_0].
\end{equation*}
Inserting this result into the second formula of \eqref{v1v2} yields
\begin{equation}\begin{aligned}&v^{(2)}\\=&u^{(2)}\otimes
q_0+hc^{(2)}\otimes q'_0+h^2 \mathcal{L}^{(2)}_{k,r}(h^2
\mathcal{L}^{(1)}_{k,r})^{-1}\otimes I_d [v^{(1)} -u^{(1)}\otimes
q_0-hc^{(1)}\otimes q'_0]\\
=&\big(u^{(2)}-A_1u^{(1)}\big)\otimes
q_0+h\big(c^{(2)}-A_1c^{(1)}\big)\otimes q'_0+A_1\otimes I_d
v^{(1)},\end{aligned}\label{v2 result}%
\end{equation}
where $A_1=\mathcal{L}^{(2)}_{k,r}(\mathcal{L}^{(1)}_{k,r})^{-1}.$
Then, by setting the matrices
$$B_1=\mathcal{L}^{(1)}_{k,r}\mathcal{P}^{(1)\intercal}_{k,r}\Omega^{(1)}_k,\ \ B_2=\mathcal{L}^{(1)}_{k,r}\mathcal{P}^{(2)\intercal}_{k,r} \Omega^{(2)}_{k},$$
substituting  \eqref{v2 result} into the first formula of
\eqref{v1v2} results in a discrete problem involving only the $r$
fundamental stages
\begin{equation*}\begin{aligned}v^{(1)}=&u^{(1)}\otimes q_0+hc^{(1)}\otimes q'_0+h^2 B_1\otimes I_d f(v^{(1)})+h^2 B_2\otimes I_d f(v^{(2)})\\
=&u^{(1)}\otimes q_0+hc^{(1)}\otimes q'_0+h^2 B_1\otimes I_d f(v^{(1)})\\
+&h^2 B_2\otimes I_d f\Big(\big(u^{(2)}-A_1u^{(1)}\big)\otimes
q_0+h\big(c^{(2)}-A_1c^{(1)}\big)\otimes q'_0+A_1\otimes I_d
v^{(1)}\Big).\\\end{aligned}\label{v1 result}%
\end{equation*}
Let
\begin{equation}\begin{aligned}\Psi(v^{(1)})=&v^{(1)}-u^{(1)}\otimes q_0-hc^{(1)}\otimes q'_0-h^2 B_1\otimes I_d f(v^{(1)})\\
&-h^2 B_2\otimes I_d f\Big(\hat{u}\otimes q_0+h\hat{c}\otimes
q'_0+A_1\otimes I_d
v^{(1)}\Big)\\\end{aligned}\label{v1 result}%
\end{equation}
with $$\hat{u}=u^{(2)}-A_1u^{(1)},\ \ \hat{c}=c^{(2)}-A_1c^{(1)}.$$
 The application of the simplified Newton method for solving
\eqref{v1 result} then gives
\begin{equation} [I_{rd}-h^2\tilde{C}\otimes J_0] \delta_{l}=-\Psi(v^{(1)}_{l}), \ \ v^{(1)}_{l+1}=v^{(1)}_{l}+\delta_{l}, \label{Newton method v1 result}%
\end{equation}
where $J_0=\frac{\partial f(q_0)}{\partial q}$ and
\begin{equation}\label{tilde Cc}\tilde{C}=B_1+B_2A_1
=\mathcal{L}^{(1)}_{k,r}\mathcal{P}^{(1)\intercal}_{k,r}\Omega^{(1)}_k+\mathcal{L}^{(1)}_{k,r}\mathcal{P}^{(2)\intercal}_{k,r}
\Omega^{(2)}_{k}\mathcal{L}^{(2)}_{k,r}(\mathcal{L}^{(1)}_{k,r})^{-1}\in\mathbb{R}^{r\times
r}.\end{equation}

The following result holds true for the matrix $\tilde{C}$.
\begin{mytheo}\label{eigenvalues  thm}
The eigenvalues of matrix $\tilde{C}$  coincide with those of matrix
$X_{r,r} $ defined in \eqref{XRR}.
\end{mytheo}
\textbf{Proof.} It follows form  \eqref{tilde Cc} that
\begin{equation*}\begin{aligned}\tilde{C}=&\mathcal{L}^{(1)}_{k,r}\mathcal{P}^{(1)\intercal}_{k,r}\Omega^{(1)}_k+\mathcal{L}^{(1)}_{k,r}\mathcal{P}^{(2)\intercal}_{k,r}
\Omega^{(2)}_{k}\mathcal{L}^{(2)}_{k,r}(\mathcal{L}^{(1)}_{k,r})^{-1}\\
=&\mathcal{L}^{(1)}_{k,r}\big(\mathcal{P}^{(1)\intercal}_{k,r}\Omega^{(1)}_k\mathcal{L}^{(1)}_{k,r}+\mathcal{P}^{(2)\intercal}_{k,r}
\Omega^{(2)}_{k}\mathcal{L}^{(2)}_{k,r}\big)(\mathcal{L}^{(1)}_{k,r})^{-1}\\
=&\mathcal{L}^{(1)}_{k,r}\big(\mathcal{P}^{\intercal}_{k,r}\Omega_k\mathcal{L}_{k,r}\big)(\mathcal{L}^{(1)}_{k,r})^{-1}\\
\sim& \mathcal{P}^{\intercal}_{k,r}\Omega_k\mathcal{L}_{k,r}=\mathcal{P}^{\intercal}_{k,r}\Omega_k\mathcal{P}_{k,r+2}\hat{X}_{r+2,r}.\end{aligned} %
\end{equation*}
By considering that $k>r$ and the quadrature formula $(c_i,b_i )$ is
exact for polynomials of degree no larger than $2k-1$, we obtain:
$$\sum\limits_{l=1}^
{k}b_l\widehat{P}_i(c_l)\widehat{P}_j(c_l)=\int_{0}^{1}\widehat{P}_i(x)\widehat{P}_j(x)dx=\delta_{ij},\
\ i=0,1,\ldots,r-1,\ j=0,1,\ldots,r+1.$$ Thus
$$\mathcal{P}^{\intercal}_{k,r}\Omega_k\mathcal{P}_{k,r+2}=\left(
    \begin{array}{ccc}
      I_r & \mathbf{0} & \mathbf{0} \\
    \end{array}
  \right)$$
and
  \begin{equation*}\begin{aligned}\tilde{C} \sim& \left(
    \begin{array}{ccc}
      I_r & \mathbf{0} & \mathbf{0} \\
    \end{array}
  \right)
 \left(
                                                                           \begin{array}{c}
                                                                             X_{r,r} \\
                                                                             \xi_{r-1}\xi_{r}e^{\intercal}_{r-1} \\
                                                                             \xi_{r}\xi_{r+1}e^{\intercal}_{r}

                                                                           \end{array}
                                                                         \right)=X_{r,r}.\end{aligned} %
\end{equation*}
\hfill  $\blacksquare$

\begin{rem} It follows from this theorem that the matrix $\tilde{C}$ has always
the same spectrum, independently of the choice of the fundamental
and silent abscissae.
  However,    the
condition number of  $\tilde{C}$ is greatly affected from this
choice. Clearly, a badly conditioned matrix $\tilde{C}$ would affect
the convergence of both the iterations \eqref{v1 result} and
\eqref{Newton method v1 result}.  Therefore, we consider a more
favorable formulation independent of the choice of the fundamental
abscissae of the discrete problem in  next subsection. \end{rem}

\subsection{Alternative formulation of the discrete problem} \label{subsec: Alternative formulation}
In order to overcome the previous drawback, the basic idea is to
reformulate the discrete problem by considering as unknowns the
coefficients. For this purpose, let us define the (block) vectors
\begin{equation} \gamma=\left(
                          \begin{array}{c}
                            \gamma_0 \\
                            \vdots \\
                            \gamma_{r-1} \\
                          \end{array}
                        \right),\ \
  \gamma_j=
\sum\limits_{l=1}^ {k}b_l\widehat{P}_j(c_l)f(v_l),\ \ j=0,1,\ldots,r-1. \label{rj}%
\end{equation}
Then the block vector $v$ in \eqref{block-matrix of methods} can be
expressed as
\begin{equation*}
v=u\otimes q_0+hc\otimes q'_0+h^2\mathcal{L}_{k,r}\otimes I_d
\gamma.
\end{equation*}
It follows from  \eqref{rj} that
\begin{equation*}\begin{aligned}\gamma=&\left(
                          \begin{array}{c}
                            \gamma_0 \\
                            \vdots \\
                            \gamma_{r-1} \\
                          \end{array}
                        \right)=\left(
                          \begin{array}{c}
                           \sum\limits_{l=1}^ {k}b_l\widehat{P}_0(c_l)f(v_l) \\
                            \vdots \\
                            \sum\limits_{l=1}^ {k}b_l\widehat{P}_{r-1}(c_l)f(v_l)  \\
                          \end{array}
                        \right)\\
=&\left(
    \begin{array}{cccc}
      b_1\widehat{P}_0(c_1) &  b_2\widehat{P}_0(c_2) & \cdots &  b_k\widehat{P}_0(c_k) \\
     b_1\widehat{P}_1(c_1) &  b_2\widehat{P}_1(c_2) & \cdots&  b_k\widehat{P}_1(c_k) \\
      \vdots & \vdots & \vdots & \vdots \\
      b_1\widehat{P}_{r-1}(c_1) &  b_2\widehat{P}_{r-1}(c_2) & \cdots &  b_k\widehat{P}_{r-1}(c_k) \\
    \end{array}
  \right)\otimes I_d\left(
           \begin{array}{c}
             f(v_1) \\
             f(v_2) \\
             \vdots \\
             f(v_k) \\
           \end{array}
         \right)\\
         =&\left(
    \begin{array}{cccc}
      \widehat{P}_0(c_1) &  \widehat{P}_0(c_2) & \cdots &  \widehat{P}_0(c_k) \\
    \widehat{P}_1(c_1) &  \widehat{P}_1(c_2) & \cdots&  \widehat{P}_1(c_k) \\
      \vdots & \vdots & \vdots & \vdots \\
      \widehat{P}_{r-1}(c_1) &  \widehat{P}_{r-1}(c_2) & \cdots &  \widehat{P}_{r-1}(c_k) \\
    \end{array}
  \right)\left(
           \begin{array}{cccc}
             b_1 &  &  & \\
              & b_2 &  &  \\
              &  & \ddots &  \\
              &  &  & b_k \\
           \end{array}
         \right)\otimes I_d
  \left(
           \begin{array}{c}
             f(v_1) \\
             f(v_2) \\
             \vdots \\
             f(v_k) \\
           \end{array}
         \right)\\
=& \mathcal{P}^{\intercal}_{k,r}\Omega_k\otimes I_d f(v).
                        \end{aligned}
                        \label{r result}%
\end{equation*}
Thus we obtain
\begin{equation*}\gamma=\mathcal{P}^{\intercal}_{k,r}\Omega_k\otimes
I_d f(u\otimes q_0+hc\otimes q'_0+h^2\mathcal{L}_{k,r}\otimes I_d
\gamma).\label{gamma}\end{equation*}

 For solving such a problem, one can still use a
fixed-point iteration
\begin{equation*}\gamma^{m+1}=\mathcal{P}^{\intercal}_{k,r}\Omega_k\otimes I_d
f(u\otimes q_0+hc\otimes q'_0+h^2\mathcal{L}_{k,r}\otimes I_d
\gamma^{m}),\ \ m=0,1,\ldots, \label{fixed-point}\end{equation*}
whose implementation is straightforward. One can also consider a
simplified-Newton iteration. Setting
\begin{equation*}\label{F gamma}
F(\gamma)\equiv\gamma-\mathcal{P}^{\intercal}_{k,r}\Omega_k\otimes
I_d f(u\otimes q_0+hc\otimes q'_0+h^2\mathcal{L}_{k,r}\otimes I_d
\gamma)=0\end{equation*} and, as before  $J_0=\frac{\partial
f(q_0)}{\partial q}$,
  it takes the form
\begin{equation} [I_{rd}-h^2\mathcal{P}^{\intercal}_{k,r}\Omega_k\mathcal{L}_{k,r}\otimes  J_0] \triangle^{m}
=-F(\gamma^{m}), \ \ \gamma^{m+1}=\gamma^{m}+ \triangle^{m}.\label{Newton method r result}%
\end{equation}
We compute that
\begin{equation*}\begin{aligned}\mathcal{P}^{\intercal}_{k,r}\Omega_k\mathcal{L}_{k,r}=\mathcal{P}^{\intercal}_{k,r}\Omega_k\mathcal{P}_{k,r+2}\hat{X}_{r+2,r}
=\left(
    \begin{array}{ccc}
      I_r & \mathbf{0} & \mathbf{0} \\
    \end{array}
  \right)
 \left(
                                                                           \begin{array}{c}
                                                                             X_{r,r} \\
                                                                             \xi_{r-1}\xi_{r}e^{\intercal}_{r-1} \\
                                                                             \xi_{r}\xi_{r+1}e^{\intercal}_{r}

                                                                           \end{array}
                                                                         \right)=X_{r,r}.\end{aligned}
\end{equation*}
Consequently, the iteration \eqref{Newton method r result} becomes
\begin{equation} [I_{rd}-h^2X_{r,r}\otimes J_0] \triangle^{m}
=-F(\gamma^{m}), \ \ \gamma^{m+1}=\gamma^{m}+ \triangle^{m}.\label{new Newton method r result}%
\end{equation}


\subsection{Blended RKN-type  Fourier collocation
methods} \label{subsec: Alternative formulation} From the arguments
in the previous subsection, it can be concluded
  that the solution of \eqref{new Newton method r result} is required at each integration step   when approximating
the   problem \eqref{common prob}. We are going to solve such an
equation by means of a blended implementation of the method. Blended
implicit methods provide a general framework for the efficient
solution of the discrete problems generated by block implicit
methods. Many researches have been done about the study of  blended
implementation of numerical methods (see, e.g.
\cite{Brugnano2002,Brugnano2004,Brugnano2007}).

In \cite{Brugnano2007}, the authors discussed the extension of
blended implicit methods for second-order problems. Following this
work, we now study the  blended RKN-type  Fourier collocation
methods for solving second-order problems \eqref{common prob}.
Following \cite{Brugnano2007},  we   consider the following linear
test equation
\begin{equation*}
q^{\prime\prime}=-\mu^2q,\ \ \mu\in \mathbb{R}.
\end{equation*}
  The iteration
\eqref{new Newton method r result}  for solving this problem leads
to a discrete problem in the form
\begin{equation}\label{1Blended  Newton method}%
(I_{rd}-\nu^2X_{r,r}\otimes I_d) \triangle
=-F(\gamma^{m})\equiv\eta_1,\ \
\nu=\mathrm{i}h\mu\equiv\mathrm{i}x,\ \ x\in \mathbb{R}.
\end{equation}
Then, we can define the following equivalent formulation of
\eqref{1Blended  Newton method}
\begin{equation}\label{2Blended  Newton method}%
\rho^2(X^{-1}_{r,r}\otimes I_d-\nu^2I_{rd}) \triangle
=\rho^2X^{-1}_{r,r}\otimes I_d\eta_1\equiv\eta_2
\end{equation}
with $\rho>0$ a free parameter. The weighting function has the form
$$\Theta(\nu)=I_r\otimes (I_d-\nu^2\rho^2I_d)^{-1}$$
and the resulting blended method corresponding to \eqref{1Blended
Newton method} is formally  given by
\begin{equation*}\label{Blended formula}%
\Theta(\nu)\eta_1+(I_{rd}-\Theta(\nu))\eta_2=[A(\nu)-\nu^2
B(\nu)]\triangle\equiv M(\nu)\triangle
\end{equation*}
with
\begin{equation*}\label{matrix Blended}%
A(\nu)=\Theta(\nu)+(I_{rd}-\Theta(\nu))\rho^2X^{-1}_{r,r}\otimes
I_d,\ \ B(\nu)=\Theta(\nu)X_{r,r}\otimes
I_d+I_{rd}-\rho^2\Theta(\nu).
\end{equation*}
 Let$$T(\triangle)=M(\nu)\triangle-\Theta(\nu)\eta_1-(I_{rd}-\Theta(\nu))\eta_2.$$
 Consequently, it induces the blended iteration
\begin{equation*}\label{blended iteration}%
\triangle_{n+1}=\triangle_{n}-\Theta(\nu)T(\triangle_{n}),
\end{equation*}
which can be rewritten as
\begin{equation*}\label{2blended iteration}%
\begin{aligned}
&N(\nu)\triangle_{n+1}=N(\nu)\triangle_{n}-T(\triangle_{n})\\
=&(N(\nu)-M(\nu))\triangle_{n}+\Theta(\nu)\eta_1+(I_{rd}-\Theta(\nu))\eta_2\equiv(N(\nu)-M(\nu))\triangle_{n}+\eta
\end{aligned}\end{equation*}
with $$N(\nu)=I_r\otimes (I_d-\nu^2\rho^2I_d)\equiv
\Theta(\nu)^{-1},\ \
\eta=\Theta(\nu)\eta_1+(I_{rd}-\Theta(\nu))\eta_2.$$

 According to the linear analysis of convergence proposed in  \cite{Brugnano2007}, the free parameter
$\rho^2$ is chosen as $$\rho^2 \equiv\min \{|\lambda| : \lambda \in
\sigma(X_{r,r})\},$$   which provides optimal convergence
properties.

A few values of $\rho^2$ are listed in Table \ref{tab1}, for the
sake of completeness.
\begin{table}$$
\begin{array}{|c|c|c|c|c|c|c|}
\hline
\text{$r$} &2  &3   &4 &5 &6   &7   \\
\hline
\text{$\rho^2$}  &6.455e-02    &3.205e-02     & 1.872e-02   &1.555e-02   &8.465e-03    &6.214e-03 \cr
 \hline
\end{array}
$$
\caption{A few values of $\rho^2$ for different $r$.} \label{tab1}
\end{table}
\subsection{Actual blended implementation} \label{sec:Actual blended
implementation} Let us now sketch the blended implementation of
RKN-type  Fourier collocation methods, when applied to a general
nonlinear system. In the case of the initial value problem
\eqref{common prob}, the previous arguments can be generalized in a
straightforward way, by considering that now the formulae
\eqref{1Blended  Newton method}-\eqref{2Blended  Newton method} and
the weighting function become
\begin{equation*}\label{blended iteration-cond}%
\begin{aligned}
&(I_{rd}-h^2X_{r,r}\otimes J_0) \triangle
=-F(\gamma^{m})\equiv\eta_1,\\
&\rho^2(X^{-1}_{r,r}\otimes I_d-h^2I_{r}\otimes J_0) \triangle
=\rho^2X^{-1}_{r,r}\otimes I_d\eta_1\equiv\eta_2,\\
&\theta=I_r\otimes(I_d-\rho^2h^2J_0)^{-1}.\end{aligned}\end{equation*}
\begin{table}$$
\begin{array}{|l|}
\hline
\text{ }     \\
\text{\ \ $\Gamma=\mathcal{P}^{\intercal}_{k,r}\Omega_k\otimes
I_d$}     \\
\text{\ \ $\Theta=h^2\mathcal{L}_{k,r}\otimes I_d$}     \\
\text{\ \ $\Upsilon=u\otimes q_0+hc\otimes q'_0$}     \\
\text{\ \ $\Lambda=\rho^2X^{-1}_{r,r}\otimes I_d$}     \\
\text{\ \ $J=\rho^2h^2J_0$}     \\
\text{\ \ $\theta=I_r\otimes(I_d-J)^{-1}$}     \\
\text{\ $\tilde{M}=\theta(I_{rd}-h^2X_{r,r}\otimes J_0)
+(I_{rd}-\theta)(\Lambda-I_{r}\otimes J)$}    \\
\text{\ \ $\gamma^0$ given\ \  \%  e.g. $\gamma^0=\Gamma f(\Upsilon)$}     \\
\text{for $l=0,1,\cdots$}    \\
\text{\ \ \ \ $v^{l}=\Upsilon+\Theta\gamma^{l}$}    \\
\text{\ \ \ \ $f^{l}=f(v^{l})$}    \\
\text{\ \ \ \ $\eta_1^{l}=-\gamma^{l}+\Gamma f^{l}$ \ \%   $-F(\gamma^{l})$}    \\
\text{\ \ \ \ $\eta_2^{l}=\Lambda\eta_1^{l}$}    \\
\text{\ \ \ \ $\triangle^{l,0}=0$}    \\
\text{\ \ \ \ \ \ \  \ for\ $j=0,1,\cdots$}    \\
\text{\ \ \ \ \ \ \ \ \ \  $\triangle^{l,j+1}=\triangle^{l,j}-\theta (\tilde{M}\triangle^{l,j}-\eta_2^{l}-\theta(\eta_1^{l}-\eta_2^{l})),$}    \\
\text{\ \ \ \ \ \ \ \ end\ $\Rightarrow$  \ returns $\triangle^{l}$}  \\
\text{\ \ \ \ $\gamma^{l+1}=\gamma^{l}+\triangle^{l}$}    \\
\text{end}    \\
  \cr
 \hline
\end{array}
$$\caption{Outer-inner iteration for the blended implementation of
RKN-type  Fourier collocation methods.} \label{tab2}
\end{table}
According to the analysis stated above, the blended iteration for
solving the first formula of \eqref{new Newton method r result} is
given as follows
\begin{equation}\label{FIN-blended iteration}%
\triangle_{n+1}=\triangle_{n}-\theta G(\triangle_{n}),\
n=0,1,\ldots,
\end{equation}
where \begin{equation*}\label{blended iteration-cond}%
\begin{aligned}
&G(\triangle)=\tilde{M}\triangle-\theta\eta_1-(I_{rd}-\theta)\eta_2,\\
&\tilde{M}=\theta(I_{rd}-h^2X_{r,r}\otimes J_0) +(I_{rd}-\theta)\rho^2(X^{-1}_{r,r}\otimes I_d-h^2I_{r}\otimes J_0).\\
\end{aligned}\end{equation*}
This iteration can be rewritten as
\begin{equation*}\label{2FIN-blended iteration}%
\begin{aligned}
&\tilde{N}\triangle_{n+1}=\tilde{N}\triangle_{n}-G(\triangle_{n})\\
=&(\tilde{N}-\tilde{M})\triangle_{n}+\theta\eta_1+(I_{rd}-\theta)\eta_2\equiv(\tilde{N}-\tilde{M})\triangle_{n}+\eta
\end{aligned}\end{equation*}
with $$\tilde{N}=I_r\otimes (I_d-\nu^2\rho^2I_d)\equiv \theta^{-1},\
\ \eta=\theta\eta_1+(I_{rd}-\theta)\eta_2.$$

From \eqref{new Newton method r result} and \eqref{FIN-blended
iteration}, we have to solve the outer-inner iteration described in
Table \ref{tab2}. A simplified (and sometimes more efficient)
procedure is given  by performing exactly 1 inner iteration (i.e.,
that with $j = 0$ in the inner cycle in Table \ref{tab2}) in the
above procedure and the corresponding  algorithm is depicted in
Table \ref{tab3}.
\begin{table}$$
\begin{array}{|l|}
\hline
\text{ }     \\
\text{\ \ $\Gamma=\mathcal{P}^{\intercal}_{k,r}\Omega_k\otimes
I_d$}     \\
\text{\ \ $\Theta=h^2\mathcal{L}_{k,r}\otimes I_d$}     \\
\text{\ \ $\Upsilon=u\otimes q_0+hc\otimes q'_0$}     \\
\text{\ \ $\Lambda=\rho^2X^{-1}_{r,r}\otimes I_d$}     \\
\text{\ \ $J=\rho^2h^2J_0$}     \\
\text{\ \ $\theta=I_r\otimes(I_d-J)^{-1}$}     \\
\text{\ $\tilde{M}=\theta(I_{rd}-h^2X_{r,r}\otimes J_0)
+(I_{rd}-\theta)(\Lambda-I_{r}\otimes J)$}    \\
\text{\ \ $\gamma^0$ given\ \  \%  e.g. $\gamma^0=\Gamma f(\Upsilon)$}      \\
\text{for $l=0,1,\cdots$}     \\
\text{\ \ \ \ $v^{l}=\Upsilon+\Theta\gamma^{l}$}    \\
\text{\ \ \ \ $f^{l}=f(v^{l})$}    \\
\text{\ \ \ \ $\eta_1^{l}=-\gamma^{l}+\Gamma f^{l}$\ \%   $-F(\gamma^{l})$}    \\
\text{\ \ \ \ $\eta_2^{l}=\Lambda\eta_1^{l}$}    \\
\text{\ \ \ \   $\triangle^{l}=\theta (\eta_2^{l}+\theta(\eta_1^{l}-\eta_2^{l})),$}    \\
\text{\ \ \ \ $\gamma^{l+1}=\gamma^{l}+\triangle^{l}$}    \\
\text{end}     \\
  \cr
 \hline
\end{array}
$$\caption{Nonlinear iteration for the blended implementation of
RKN-type  Fourier collocation methods.} \label{tab3}
\end{table}
\section{Numerical tests} \label{sec:Numerical tests}
In this section, a couple of numerical examples are  shown   to put
into evidence on the features and effectiveness of the methods.
 As an example of the RKN-type  Fourier collocation
methods, we choose $k=4$ and the following  Gauss--Legendre's
 quadrature
\begin{equation*}\begin{aligned}
&c_1=\frac{1+\sqrt{\frac{3}{7}+\frac{2\sqrt{\frac{6}{5}}}{7}}}{2},\
\
c_2=\frac{1+\sqrt{\frac{3}{7}-\frac{2\sqrt{\frac{6}{5}}}{7}}}{2},\\
&c_3=\frac{1-\sqrt{\frac{3}{7}-\frac{2\sqrt{\frac{6}{5}}}{7}}}{2},\
\
c_4=\frac{1-\sqrt{\frac{3}{7}+\frac{2\sqrt{\frac{6}{5}}}{7}}}{2}, \\
&b_1=\frac{1}{2}(\frac{1}{2}-\frac{1}{6}\sqrt{\frac{5}{6}}),\ \ \ \ \ b_2=\frac{1}{2}(\frac{1}{2}+\frac{1}{6}\sqrt{\frac{5}{6}}),\\
&b_3=\frac{1}{2}(\frac{1}{2}+\frac{1}{6}\sqrt{\frac{5}{6}}),\ \ \ \ \ b_4=\frac{1}{2}(\frac{1}{2}-\frac{1}{6}\sqrt{\frac{5}{6}}).\\
\label{GL}
\end{aligned}\end{equation*}
Then $r=2$ is chosen in \eqref{methods0}, and the corresponding
method is denoted by RKN-TFC. According to the analysis of
\cite{wang-2016}, this method is of order four. In order to show the
efficiency and robustness of the fourth-order method, the
integrators we select for comparison are:

\begin{itemize}\itemsep=-0.2mm

\item RKN-TFC-B: the RKN-type Fourier collocation method RKN-TFC of
order four using the blended iteration described in Table
\ref{tab3};

\item RKN-TFC-F: the   method RKN-TFC
using the fixed-point iteration;

\item DIRKN-F: the A-stable   diagonally
implicit RKN method DIRKN$(2)4_{12-22}$ of order four in
\cite{franco2009} using the fixed-point iteration.

\end{itemize}

 In the practical
computations, we set $10^{-16}$ as the error tolerance and $10^4$ as
the maximum number of each iteration.
\begin{table}
\newcommand{\tabincell}[2]{\begin{tabular}{@{}#1@{}}#2\end{tabular}}
  \centering
  \begin{tabular}{|c|c|c|c|c|c|c|}\hline
Methods $(t,h)$ & CPU time &Iterations& \tabincell{c}{Solution\\
error}&\tabincell{c}{Hamiltonian\\ error}&\tabincell{c}{Angular\\
momentum\\ error}
\\\hline  RKN-FCM-B $(50,0.4)$      &0.029 &1423&$-$2.149&$-$9.248 &  $-$9.069\\
  RKN-FCM-F $(50,0.4)$     & 0.709&71438 &0.300&$-$3.143& $-$4.126\\
    DIRKN-F $(50,0.4)$     &  0.132& 9676 &$-$0.544 &$-$0.627 & $-$0.753 \\
\hline  RKN-FCM-B  $(50,0.2)$    &0.064 &3028&$-$3.354&$-$11.700 & $-$11.524\\
 RKN-FCM-F  $(50,0.2)$      &0.531 &52204 &$-$0.196&$-$4.930&$-$5.558\\
     DIRKN-F $(50,0.2)$     &0.111 &  5966 &$-$1.202 &$-$3.056 & $-$3.057 \\
\hline  RKN-FCM-B  $(50,0.1)$     & 0.080  &3285&$-$4.558&$-$14.002  & $-$13.875\\
 RKN-FCM-F  $(50,0.1)$   & 0.257 &23490&$-$0.711&$-$6.432& $-$7.038\\
     DIRKN-F $(50,0.1)$     & 0.085 &  5486 & $-$3.006&$-$4.887 & $-$4.88\\
\hline   RKN-FCM-B  $(100,0.4)$       &  0.090 &3841&$-$1.879&$-$8.658&  $-$8.479\\
  RKN-FCM-F  $(100,0.4)$    & 1.692 &162856  &$-$0.255& $-$2.785& $-$3.822\\
      DIRKN-F $(100,0.4)$     & 0.365 & 26334 &$-$0.297 &$-$0.627& $-$0.753 \\
\hline  RKN-FCM-B  $(100,0.2)$    &0.137  &7048&$-$3.085&$-$11.109& $-$10.932\\
  RKN-FCM-F  $(100,0.2)$ & 0.866  &84425 &0.129&$-$4.413&$-$5.255\\
      DIRKN-F $(100,0.2)$     & 0.134& 9952  &$-$0.603 &$-$2.750 & $-$2.752\\
\hline  RKN-FCM-B $(100,0.1)$     &0.177  &7573&$-$4.289&$-$13.461& $-$13.331\\
 RKN-FCM-F  $(100,0.1)$   & 0.529 &46986&$-$0.484&$-$5.900&$-$6.736\\
     DIRKN-F $(100,0.1)$     & 0.164 & 9981  &$-$2.441 &$-$4.586 & $-$4.587\\
\hline
\end{tabular}
 \caption{Results for Problem 1.} \label{tab 4}
\end{table}

 \vspace{1mm}\noindent\textbf{Problem 1.} Consider  the following perturbed
Kepler's problem
\begin{equation*}
\begin{aligned}
&q_1^{\prime\prime}=-\dfrac{q_1}{(q_1^2+q_2^2)^{3/2}}-\dfrac{(2\epsilon+\epsilon^2)q_1}{(q_1^2+q_2^2)^{5/2}},
 \quad q_1(0)=1,\quad q_1'(0)=0,\\
&q_2^{\prime\prime}=-\dfrac{q_2}{(q_1^2+q_2^2)^{3/2}}-\dfrac{(2\epsilon+\epsilon^2)q_2}{(q_1^2+q_2^2)^{5/2}},
 \quad q_2(0)=0,\quad q_2'(0)=1+\epsilon.
\end{aligned}
\end{equation*}
The exact solution   is $$q_1(t)=\cos(t+\epsilon t),\
q_2(t)=\sin(t+\epsilon t).$$ The Hamiltonian is
$$H=\dfrac{1}{2}(q_1'^2+q_2'^2)-\dfrac{1}{\sqrt{q_1^2+q_2^2}}-\dfrac{(2\epsilon+\epsilon^2)}{3(q_1^2+q_2^2)^{3/2}}.$$
The system also has the angular momentum $L=q_1q_2'-q_2q_1'$ as a
first integral. We take the parameter value $\epsilon=10^{-3}$.

We   solve the problem in the intervals $[0, 50]$ and  $[0, 100]$
 with different stepsizes  $h=0.4,0.2,0.1.$    Table \ref{tab 4}
lists  the CPU time, the total numbers of iterations, the solution
errors, Hamiltonian errors and
 angular momentum  errors.  It is noted here that all the errors
are given by the logarithm of the corresponding results.

 \vspace{1mm} \noindent\textbf{Problem 2.}  The H\'{e}non-Heiles Model is created for
describing stellar motion and it can be expressed by the following
form
\[
\begin{aligned}
&q''_1(t)=-q_1(t)-2q_1(t)q_2(t),\ \ \ q_1(0)=\sqrt{\dfrac{11}{96}},\ q'_1(0)=0,\\
&q''_2(t)=-q_2(t)-q_1^2(t)+q_2^2(t),\ q_2(0)=0,\ \ \ \ \ \   q'_2(0)=\dfrac{1}{4}.\\
\end{aligned}
\]
 The Hamiltonian function of the system is given by
\[
H(p,q)=\dfrac{1}{2}(q_1'^2+q_2'^2)+\dfrac{1}{2}(q_{1}^{2}+q_{2}%
^{2})+q_{1}^{2}q_{2}-\dfrac{1}{3}q_{2}^{3}.
\]

\begin{table}
\newcommand{\tabincell}[2]{\begin{tabular}{@{}#1@{}}#2\end{tabular}}
  \centering
  \begin{tabular}{|c|c|c|c|c|c|}\hline
Methods $(t,h)$ & CPU time &Iterations& \tabincell{c}{Solution\\
error}&\tabincell{c}{Hamiltonian\\ error}
\\\hline  RKN-FCM-B $(50,0.1)$      &0.063 &2989& $-$5.806&$-$8.915 \\
  RKN-FCM-F $(50,0.1)$     & 0.048 &3004 &$-$2.145& $-$5.170\\
    DIRKN-F $(50,0.1)$     &  0.058&    3979  &$-$5.545&$-$5.547  \\
\hline  RKN-FCM-B  $(50,0.05)$    & 0.113&4996&$-$7.010&$-$10.121 \\
 RKN-FCM-F  $(50,0.05)$      & 0.085&5000 & $-$2.754&$-$6.005\\
     DIRKN-F $(50,0.05)$     &0.092 & 6187 &$-$6.653 &$-$7.054  \\
\hline  RKN-FCM-B  $(50,0.025)$     & 0.198 &8012&$-$8.214&$-$11.325 \\
 RKN-FCM-F  $(50,0.025)$   & 0.158 &8906&$-$3.359& $-$6.711\\
     DIRKN-F $(50,0.025)$     & 0.174 &   10018 &  $-$7.828& $-$8.560 \\
\hline   RKN-FCM-B  $(100,0.1)$       & 0.108 &5981&$-$5.301& $-$7.900\\
  RKN-FCM-F  $(100,0.1)$    &0.080 &6004   & $-$1.654&$-$4.338\\
      DIRKN-F $(100,0.1)$     & 0.091 & 7958 & $-$4.247 &$-$5.245 \\
\hline  RKN-FCM-B  $(100,0.05)$    &0.206& 9996&$-$6.504&$-$9.105 \\
  RKN-FCM-F  $(100,0.05)$ & 0.162 &10000 & $-$2.253&$-$4.927 \\
      DIRKN-F $(100,0.05)$     &  0.167&  12366  &$-$5.747 &$-$6.750  \\
\hline RKN-FCM-B $(100,0.025)$     &0.368 &16025& $-$7.708&$-$10.309  \\
 RKN-FCM-F  $(100,0.025)$   & 0.300  &17801& $-$2.854& $-$5.528 \\
     DIRKN-F $(100,0.025)$     &  0.318 &  20035  & $-$7.193 &$-$8.253 \\
\hline
\end{tabular}
\caption{Results for Problem 2.} \label{pro2 tab}
\end{table}

We    solve the problem in the intervals $[0, 50]$ and $[0, 100]$
with different stepsizes  $h=0.4,0.2,0.1$.    See Table \ref{pro2
tab} for  the CPU time, the total numbers of iterations, the
solution errors and Hamiltonian errors.

From the numerical tests, one can then conclude that the proposed
blended implementation of  RKN-type Fourier collocation methods
turns out to be more robust and efficient than the same method and
the RKN method using the fixed-point iteration.

 \section{Conclusions} \label{sec:conclusions}
In this paper we propose and analyze an efficient iterative
procedure for solving the second-order differential equations
\eqref{common prob} generated by the application of RKN-type Fourier
collocation methods. The proposed implementation turns out to be
robust and efficient. Two numerical tests confirm the effectiveness
of the proposed iteration when numerically solving second-order
differential equations.

Last but not least, it is noted that there are still some issues
which will be further considered.
\begin{itemize}
\item  A new splitting procedure has been developed recently  for the implementation of implicit methods
and we refer the reader to \cite{Brugnano2013,Brugnano2015} for
example. We will consider the novel technique for RKN-type Fourier
collocation methods  in a future research.

\item Trigonometric Fourier collocation methods (TFCMs) have been formulated in \cite{wang-2016} for
solving $q^{\prime\prime}+Mq=f(q)$ with a matrix $M$ and the
coefficients of the methods depend on  matrix-valued functions.
Therefore, the implementation of these methods would be different
and complicated. Another issue for future exploration is the
research of efficient implementation for the TFCMs.

\item The shifted Legendre polynomials are chosen   in \cite{wang-2016} as  an orthonormal
basis to give the formulation of the methods. It is noted that a
different choice of the orthonormal basis would modify the scheme of
the methods  as well as the analysis of the implementation, which
will be considered in future investigations.
\end{itemize}

\end{document}